\documentclass[12pt]{article}

\usepackage{amsmath,amsthm}
\usepackage{amssymb}

\newtheorem{thm}{\small\bf Theorem}
\newtheorem{lem}{\small\bf Lemma}

\numberwithin{equation}{section}

\newcommand{\be}{\begin{eqnarray}}
\newcommand{\ee}{\end{eqnarray}}
\newcommand{\bd}{\begin{displaymath}}
\newcommand{\ed}{\end{displaymath}}

\newcommand{\R}{\mathbb{R}}

\title{A note on Blasius type boundary value problems}
\author{Grzegorz Andrzejczak \\Magdalena Nockowska-Rosiak \\ Bogdan Przeradzki  \\ Technical University of Lodz, Poland}

\begin{document}
\maketitle
\begin{abstract}
The existence and uniqueness of a solution to a generalized Blasius equation with asymptotic boundary conditions
are proved. A new numerical approximation method is proposed.
\end{abstract}
\textit{keywords:} Blasius equation, shooting method
\\\textit{2010 MSC:} 34B15, 34B40, 34D05, 65L10 
\section{Introduction}

We study the BVP of the form:
\be            \label{equ}
x'''+cx^p\cdot x''=0,\qquad x(0)=0=x'(0),\quad \lim_{t\to\infty}x'(t)=\beta,
\ee
where $p\ge 1,$
$c$ and $\beta$ are positive constants. The problem is motivated by the classical Blasius equation describing the velocity profile of the fluid in the boundary layer where $c=\frac{1}{2},$ $p=\beta =1$. The Blasius equation is a basic equation in fluid mechanics which appears in the study of the flow of an incompressible viscous fluid over semi infinite plane. Blasius (\cite{BL}) used a similarity transform technique to convert the partial differential equation into his famous ordinary differential equation
\be     \label{equ-Bl}
x'''+\tfrac{1}{2}x\cdot x''=0,\qquad x(0)=0=x'(0),\quad \lim_{t\to\infty}x'(t)=1,
\ee
where $x$ is the stream function $x=\frac{\Psi}{\sqrt{2U\nu y_1}}$, $U$ is the fluid velocity, $\nu$ is the fluid kinematic viscosity and $t$ is the similarity variable defined as $t=y_2\sqrt{\frac{U}{2\nu y_1}}$, where $y_1,y_2$ are Cartesian coordinates with $y_1$ pointing along the free stream direction and $y_2$ perpendicular to $y_1$. We refer to \cite{Bo1,Bo2} for an excellent introduction to the problem.
A series expansions method was used to solve \eqref{equ-Bl} by Blasius. There has been appeared many analytical and numerical methods handling this problem since the Blasius's work, \cite{PT,Ran} for instance.
\par
In the first part of this paper, the existence and uniqueness of \eqref{equ} will be analytically proved by changing the boundary value problem to an initial problem. Using the obtained estimates we will be able to find the value of $a=x''(0)$ which guarantees that the solution $x_a$ on an initial problem
\be     \label{equ-ini}
x'''+cx^p\cdot x''=0,\qquad x(0)=0=x'(0),\quad x''(0)=a,
\ee
is the solution of \eqref{equ} we are looking for. In the second part of the article, a new numerical approximation method is proposed.

\section{Auxiliary lemmas}

Let $x_a$ stand for the unique solution satisfying initial conditions $$x(0)=0=x'(0),\qquad x''(0)=a.$$
If $a<0,$ then $x_a$ and $x_a''$ are  negative for small $t$'s thus the solution is concave and negative for all arguments
and it cannot solve (\ref{equ}). For $a=0$ we have a trivial solution $x_a\equiv 0$ and the seeking solution can be obtained for $a>0.$

\begin{lem}
The $x_a$ is defined for all $t\ge0.$
\end{lem}
{\sl Proof.} $x_a''$ cannot vanish at any point $t_0$ by the uniqueness of solutions of initial value problems: $x(t)=c_1t+c_2$ solves our ODE.
Hence dividing the equation by $x_a''$ and integrating on $[0,t]$ we have
\be \label{2der}
x_a''(t)=a\exp\left( -c\int_0^t x_a(s)^p\, ds \right),
\ee
which implies
\be   \label{1der}
x_a'(t)=a\int_0^t \exp \left( -c\int_0^s x_a(\tau )^p\, d\tau \right)\, ds.
\ee
Integrating once more and applying the Fubini Theorem we get
\be \label {0der}
x_a(t)=a\int_0^t (t-s)\exp \left( -c\int_0^s x_a(\tau )^p\, d\tau \right)\, ds.
\ee
By (\ref{2der}),(\ref{1der}) and (\ref{0der}) we have apriori estimates:
$$0<x_a(t)<\frac{1}{2}at^2, \qquad 0<x_a'(t)<at, \qquad 0<x_a''(t)<a$$
for any $t>0.$ It follows (\cite{PSV}, p. 146) that $x_a$ is extendable to $[0,\infty).$

\vspace{.3cm}

\begin{lem} \label{last}
For any $a>0,$ $\lim_{t\to\infty} x_a''(t)=0$ and there exists a finite and positive limit $h(a):=\lim_{t\to\infty} x_a'(t).$
\end{lem}
{\sl Proof.} Since $x_a''>0$ and $x_a'>0$, then $\lim_{t\to+\infty}x_a(t)=+\infty$. Moreover, since $x''_a>0$, $x_a>0$ and $x'''_a=-cx_a^px''_a$, then $x''_a$ is a decreasing function so $\lim_{t\to\infty} x_a''(t)=g_a\in[0,a)$. Suppose $g_a>0$. From $x'''_a=-cx_a^px''_a$ we get that
\be\label{pom}
\lim_{t\to+\infty}x'''_a(t)=-\infty.
\ee
On the other hand
$$
\forall \, t\geq0\ \exists\, s_t\in(t,t+1)\quad x''_a(t+1)-x''_a(t)=x'''_a(s_t)
$$
and hence $\lim_{t\to+\infty}x'''_a(t)=0$, which contradicts \eqref{pom}.
\par
Since $x_a''>0,$ then $x_a'$ is an increasing function (so the limit defining $h(a)$ exists, possibly infinite). From $\lim_{t\to+\infty}x_a(t)=+\infty$ we get there exists $t_a>0$ such that $cx_a(t)^p>1$ for $t>t_a.$ From (\ref{1der}), we obtain for $t>t_a,$
$$x_a'(t)\le a\int_0^t \exp\left( -c\int_0^{t_a} x_a(\tau)^p\, d\tau\right) \exp\left(-\int_{t_a}^s \, d\tau\right)\, ds\le$$
$$\le a\int_0^{t_a}\, ds +a \int_{t_a}^t  \textrm{e}^{-(s-t_a)}\, ds \le at_a+  a\textrm{e}^{ t_a} \int_{t_a}^{\infty}  \textrm{e}^{-s}\, ds =a(t_a+1).$$
Thus
\be \label{h-up}
h(a)\le a(t_a+1)
\ee
 for any $a>0.$

\vspace{.3cm}

\begin{lem}
For any $a>0,$ there exists a finite and positive limit \newline $\mu(a):=\lim_{t\to\infty} (h(a)t-x_a(t)).$
It means that the graph of $x_a$ has a slant asymptote and the following estimates hold:
\be \label{est}
\max (0,h(a)t-\mu(a))\le x_a(t)\le h(a)t.
\ee
\end{lem}
{\sl Proof.} The function $t\mapsto h(a)t-x_a(t)$ is increasing, hence the limit from the assertion exists but it can be infinite. Suppose it equals $+\infty.$
By the arguments from the proof of the  lemma \ref{last}, we have $x_a'''(t)\le  -x_a''(t)$ for $t\ge t_a.$ Integrating this inequality from $s$ to $+\infty$
and using the fact $x_a''(+\infty)=0,$ we get
$$-x_a''(s)\le -h(a)+x_a'(s).$$
Next integration from $t_a$ to $t$ leads to the following inequality
$$x_a'(t_a)-x_a'(t)\le -h(a)(t-t_a)+x_a(t)-x_a(t_a)$$
or equivalently
$$h(a)t-x_a(t)\le h(a)t_a-x_a(t_a)+h(a)-x_a'(t_a).$$
Thus
\be  \label{mu}
0< \mu (a)\le h(a)t_a-x_a(t_a)+h(a)-x_a'(t_a).
\ee
The last part of the assertion is a simple consequence.

\vspace{.3cm}
For an upper bound on $h(a)$, $\mu(a)$ depending explicitily on $a$ we use $x_a(t)\le at^2/2$ to (\ref{0der}). Hence,
$$x_a(t)\ge a\int_0^t (t-s)\exp \left( -\int_0^s c\frac{a^p}{2^p}\tau ^{2p}\, d\tau \right)\, ds= a\int_0^t (t-s)\exp \left( -c\frac{a^p}{(2p+1)2^p} s^{2p+1}\right)\, ds.$$
One can easily show that the function $$\varphi(t):=\int_0^t (t-s)\exp\left( -ks^{\alpha}\right)\, ds, \ k=\frac{ca^p}{2^p(2p+1)}, \ \alpha=2p+1$$
has a similar behaviour as $x_a$ in the sense that its graph has an asymptote $x=h^{\ast} t-\mu^{\ast},$
where
\be\label{exp}
h^{\ast}=\int_0^{\infty}\exp\left( -ks^{\alpha}\right)\, ds=\frac{\Gamma (1/\alpha)}{\alpha k^{1/{\alpha}}},\qquad \mu^{\ast}=\int_0^{\infty} s\exp\left( -ks^{\alpha}\right)\, ds=\frac{\Gamma (2/\alpha)}{\alpha k^{2/{\alpha}}}
\ee
and its graph sits above this line. Hence,
$$x_a(\tau)\ge a^{1-\frac{p}{2p+1}}\cdot \Gamma(\tfrac{1}{2p+1})(\tfrac{2^p}{c\cdot(2p+1)^{2p}})^{\frac{1}{2p+1}}\cdot\tau -a^{1-\frac{2p}{2p+1}}\cdot \Gamma(\tfrac{2}{2p+1})(2p+1)^{\frac{1-2p}{2p+1}}(\tfrac{2^p}{c})^{\frac{2}{2p+1}}.$$
\be\label{osz-1}
x_a(\tau)\ge c_2\cdot a^{\frac{p+1}{2p+1}}\cdot\tau -c_3\cdot a^{\frac{1}{2p+1}},
\ee
where
\be\label{constant}
c_2:=\Gamma(\tfrac{1}{2p+1})(\tfrac{2^p}{c\cdot(2p+1)^{2p}})^{\frac{1}{2p+1}}, \quad c_3:=\Gamma(\tfrac{2}{2p+1})(2p+1)^{\frac{1-2p}{2p+1}}(\tfrac{2^p}{c})^{\frac{2}{2p+1}}.
\ee
Now, we are able to get appropriate estimates for $h(a).$
\begin{lem}  \label{down}
For any $a>0,$
\be \label{h}
c_2\cdot a^{(p+1)/(2p+1)}\le h(a) \le c_1\cdot a^{(p+1)/(2p+1)} ,
\ee
where
$$
c_1:=\frac{c_3}{c_2}+\frac{\Gamma(1/(p+1))}{c^{1/(p+1)}\cdot (c_2(p+1))^{p/(p+1)}}.$$
\end{lem}
{\sl Proof.} For a lower bound we apply the estimate $x_a(\tau)\le \frac{1}{2}a\tau ^2$ to the equality
\be \label{wzor-h}
h(a)=a\int_0^{\infty} \exp \left( -c\int_0^s x_a(\tau )^p\, d\tau \right)\, ds.
\ee
This leads to the inequality
$$h(a)\ge a\int_0^{\infty} \exp\left( -\frac{ca^p}{2^p(2p+1)} s^{2p+1} \right)\, ds.$$
\eqref{exp} for $k=\frac{ca^p}{2^p(2p+1)}$, $\alpha =2p+1$ gives the lower bound on $h$.
\\For an upper bound on $h$ we use the lower estimate of $x_a$ - \eqref{osz-1} to the equality \eqref{wzor-h} and we get
\begin{align*}
&h(a)\le a \int_0^{c_3/c_2\cdot a^{-p/(2p+1)}} \, ds +
\\
&a\int_{c_3/c_2\cdot a^{-p/(2p+1)}}^{\infty} \exp\left( -c\int_{c_3/c_2\cdot a^{-p/(2p+1)}}^s \left( c_2a^{\frac{p+1}{2p+1}}\tau -c_3 a^{\frac{1}{2p+1}}\right) ^p \, d\tau\right)\, ds
\\
&\le \frac{c_3}{c_2}a^{\frac{p+1}{2p+1}} + \frac{1}{c_2} a^{\frac{p}{2p+1}}\int_0^{\infty} \exp\left( -\frac{c}{(p+1)c_2}a^{-\frac{p+1}{2p+1}}t^{p+1}\right)\, dt,
\end{align*}
where we used linear substitutions twice. At last, using \eqref{exp} for $k=\frac{c}{(p+1)c_2}a^{-\frac{p+1}{2p+1}}$, $\alpha =p+1$  we have for any $a>0,$
$$h(a)\le \bigg(\frac{c_3}{c_2} + \frac{\Gamma(1/(p+1))}{c^{1/(p+1)}\cdot (c_2(p+1))^{p/(p+1)}}\bigg)\cdot a^{\frac{p+1}{2p+1}}.$$
The next lemma presents estimates for $\mu(a)$.
\begin{lem}
For any $a>0,$ constant $\mu (a)$ satisfies the following estimates:
$$c_4\cdot a^{1/(2p+1)}\le \mu(a) \le c_5\cdot a^{1/(2p+1)},$$
where
\begin{align*}
&c_4:= \frac{2^{2p/(2p+1)}\Gamma(2/(2p+1))}{(2p+1)^{(2p-1)/(2p+1)}c^{2/(2p+1)}},
\\ &c_5:=\tfrac{1}{c^2_2}\big[\tfrac{c^2_3}{2}+(\tfrac{c_2}{c})^{2/(p+1)}(p+1)^{(1-p)/(1+p)}\Gamma(\tfrac{2}{p+1})+c_3(\tfrac{c_2}{c\cdot(p+1)^p})^{1/(p+1)}\Gamma(\tfrac{1}{p+1})\big].
\end{align*}
\end{lem}
{\sl Proof.} From (\ref{0der}) and (\ref{wzor-h}) we get
\be \label{wzor-mu}
\mu (a)=a\int_0^{\infty} s\exp\left( -c\int_0^s x_a(\tau)^p\, d\tau \right)\, ds.
\ee
Using the estimate $x_a(\tau)\le\tfrac{1}{2}a\tau^2$ we get
$$\mu(a)\ge a\int_0^{\infty}s\cdot \exp\left( -\frac{ca^p}{2^p(2p+1)} s^{2p+1}\right)\, ds.$$
Using \eqref{exp} for $k=\frac{ca^p}{2^p(2p+1)}$, $\alpha =2p+1$ we obtain the lower bound.
On the other hand, from \eqref{osz-1}, \eqref{wzor-mu}  we have
\begin{align*}
&\mu(a)\le a \int_0^{c_3/c_2\cdot a^{-p/(2p+1)}} s\, ds +
\\
&a\int_{c_3/c_2\cdot a^{-p/(2p+1)}}^{\infty} s\exp\left( -c\int_{c_3/c_2\cdot a^{-p/(2p+1)}}^s \left( c_2a^{\frac{p+1}{2p+1}}\tau -c_3 a^{\frac{1}{2p+1}}\right) ^p \, d\tau\right)\, ds\le
\\
&\le \frac{c_3^2}{2c_2^2}a^{\frac{1}{2p+1}} + \frac{1}{c_2^2}a^{\frac{-1}{2p+1}} \int_0^{\infty} t\exp\left( -\frac{c}{(p+1)c_2}a^{-\frac{p+1}{2p+1}}t^{p+1}\right)\, dt+
\\
&+\frac{c_3}{c_2^2} \int_0^{\infty} \exp\left( -\frac{c}{(p+1)c_2}a^{-\frac{p+1}{2p+1}}t^{p+1}\right)\, dt,
\end{align*}
where we used linear substitutions twice. At last, using \eqref{exp} for $k=\frac{c}{(p+1)c_2}a^{-\frac{p+1}{2p+1}}$, $\alpha =p+1$  we have for any $a>0$ we get the upper bound.


\section{Main results}

Now, we are able to prove the existence of a solution to (\ref{equ}).
\begin{thm}
The BVP (\ref{equ}) has a solution for any $\beta \ge 0.$
\end{thm}
{\sl Proof.} The function $h:[0,\infty)\to\mathbb{R}$ is continuous on $(0,+\infty)$ by the continuous dependence of
solutions of ODEs on initial conditions and locally uniform convergence of the integral
$$\int_0^{\infty} \exp\left( -c\int_0^s x_a(\tau)^p\, d\tau \right)\, ds.$$
By the estimates from lemma \ref{down}, we have
$$\lim_{a\to 0^+} h(a)=0,\qquad \lim_{a\to\infty} h(a)=+\infty.$$
Thus, for any $\beta>0,$ there exists $a>0$ such that $h(a)=\beta.$ For $\beta=0,$ it is obvious.

\vspace{.3cm}

Finally, the uniqueness of the solution of \eqref{equ} will be proved by using the ideas from \cite{BH}. For any $a>0$ consider the one-to-one function $v_a:[0,h(a)^2)\to [0,\infty)$ such that $v_a(x_a'(t)^2)=x_a(t)$ for each $t\ge 0.$ It is well defined since $x_a$ and $x_a'$ are increasing functions and it belongs to $C^2(0,h(a)^2).$ Substituting $y=x'(t)^2,$ we shall find an ODE satisfied by $v$ (we omit subscript $a$ for simplicity).
$$x(t)=v(y),\qquad x'(t)=v'(y)2x'(t)x''(t)$$
hence, $$x''(t)=\frac{1}{2v'(y)}, \qquad x'''(t)=-\frac{v''(y)}{2v'(y)^2}\cdot 2x'(t)x''(t)=-\frac{v''(y)\sqrt{y}}{2v'(y)^3}.$$
Put $x$ and its derivatives in our ODE and find
\be \label{v}
v''(y)=\frac{cv(y)^p v'(y)^2}{\sqrt{y}}.
\ee
From boundary conditions on $x$ we get
\be \label{bound-v}
v(0)=0,\qquad v'(0)=\frac{1}{2a},\qquad \lim_{y\to h(a)^2 -} v(y)=+\infty.
\ee
Now, we are in position to prove
\begin{thm}
The solution of (\ref{equ}) is unique.
\end{thm}
{\sl Proof.} We need to show that the function $h$ is one-to-one. Suppose that $h(a_1)=h(a_2),$ $a_2>a_1$ and take $v_1$ and $v_2$ obtained by $x_{a_1}$ and $x_{a_2},$ respectively, that is $v_i$ satisfies (\ref{v}) with boundary conditions (\ref{bound-v}) (for $a=a_i,$ $i=1,2).$ Put $w=v_1-v_2.$ Then $w(0)=0,$ $w'(0)=\frac{a_2 -a_1}{2a_1 a_2}>0.$ Notice that $w'>0$ on the whole interval $(0,h(a_1)^2)$ -- both function are defined on the same interval.

In fact, if it is not true, then there exists $s$ in this interval such that $w'>0$ on $(0,s)$ and $w'(s)=0.$ Hence $w(s)>w(0)=0$ and
$$w''(s)=\lim_{\xi\to 0+} \frac{w'(s-\xi)-w'(s)}{-\xi}\le 0.$$
On the other hand, $$w''(s)=v_1''(s)-v_2''(s)=cs^{-1/2}(v_1(s)^p-v_2(s)^p)v_{a_i}'(s)^2$$
since $w'(s)=0$ implies $v_{a_1}'(s)=v_{a_2}'(s).$ But $v_1(s)^p>v_2(s)^p$ from $w(s)>0$ and this gives $w''(s)>0$ -- a contradiction. Thus, we have $w>0$ and $w'>0$ on $(0,h(a_1)^2).$

Set $V_i=1/v_i',$ $i=1,2$ and $W=V_1-V_2.$ We have, for any $y\in (0,h(a_1)^2),$
$$W'(y)=V_1'(y)-V_2'(y)=-\frac{v_1''(y)}{v_1'(y)^2}+\frac{v_2''(y)}{v_2'(y)^2}=c\frac{v_2(y)^p-v_1(y)^p}{\sqrt{y}}<0$$
from $w(y)>0.$ Hence,
$$W(y)<W(0)=\frac{1}{v_1'(0)}-\frac{1}{v_2'(0)}=2(a_1-a_2)$$
and $$\lim_{y\to h(a_1)^2 -}W(y)\le 2(a_1-a_2)<0.$$
On the other hand, $V_i(x_{a_i}'(t)^2)=2x_{a_i}''(t)$ implies
$$\lim_{y\to h(a_1)^2 -}V_i(y)=\lim_{t\to\infty} V_i(x_{a_i}'(t)^2)=2\lim_{t\to\infty} x_{a_i}''(t)=0$$
and, therefore, $$\lim_{y\to h(a_1)^2 -}W(y)=0$$
which contradicts the previous inequality.

\section{Numerical approach}

All numerical methods cannot work on the infinite interval $[0,\infty)$ and we do not know the exact value of $a=x''(0)$ for the solution. Our earlier results
make possible to find a finite interval $[0,T]$ for any positive value $\epsilon$ of the error control tolerance such that
$$x''(T)<\epsilon,\qquad h(a) -x'(T)<\epsilon,\qquad x(T)-(h(a) T-\mu (a))<\epsilon.$$
Since all these functions decrease, all three inequalities hold for any $t>T.$  First, by using estimates (\ref{h}), we can find an interval $[a_{min},a_{max}]$
such that $h(a_{min})<\beta <h(a_{max}).$ Next, by (\ref{2der}), we need
$$a_{max} \exp \left( -c\int_0^T x_{a_{min}}(\tau)^p\, d\tau \right)<\epsilon,$$
by (\ref{wzor-h}), we should have
$$a_{max} \int_T^{\infty} \exp\left( -c\int_0^s x_{a_{min}}(\tau)^p\, d\tau\right)ds <\epsilon,$$
and by(\ref{0der}) and (\ref{wzor-mu}), we get
$$-a_{max}\int_T^{\infty} (T-s)\exp\left( -c\int_0^s x_{a_{min}}(\tau)^p\, d\tau\right)ds <\epsilon.$$
We do not know the function $x_{a_{min}}$ but we can use estimate (\ref{osz-1}) to get
\be\label{2-der}
a_{max} \exp \left( -c\int_0^T  \left(c_2\cdot a_{min}^{\frac{p+1}{2p+1}}\cdot \tau -c_3\cdot a_{min}^{\frac{1}{2p+1}}\right)^p\, d\tau \right)<\epsilon,
\ee
\be \label{1-der}
a_{max} \int_T^{\infty} \exp\left( -c\int_0^s \left(c_2\cdot a_{min}^{\frac{p+1}{2p+1}}\cdot \tau -c_3\cdot a_{min}^{\frac{1}{2p+1}}\right)^p\, d\tau\right) ds <\epsilon,
\ee
\be \label{0-der}
-a_{max}\int_T^{\infty} (T-s)\exp\left( -c\int_0^s \left(c_2\cdot a_{min}^{\frac{p+1}{2p+1}}\cdot \tau -c_3\cdot a_{min}^{\frac{1}{2p+1}}\right)^p\, d\tau\right)ds <\epsilon.
\ee

We start with a family of initial value problems
\be \label{numer}
x_a(t)\;\; \textrm{for}\;\; t \in [0,T],\qquad x_a(0) = 0 = x_a'(0),\qquad x''_a(0) = a,\;\; a \in [a_{min}, a_{max}].
\ee
If we approximate this solution in $[0,T]$ with an error less than $\epsilon,$ then the best approximation of $x$ in $[T,\infty )$ is
$$x(t)=\beta t+\left( x_a(T)-\beta T\right).$$
The lower and upper bounds for the second derivative describe the shooting window -
for each $a$ the only solution in this direction exists at $t=T,$ and the computed value of $x'(T)$ is more and more close to the expected limit value $\beta.$ As long as $\beta$ is contained between the computed values $x'(T)$ of the best two shots, we apply the classical bisection method:

If $y$ and $z$ are solutions such that  $y''(0)<z''(0),$ and there is $y'(T)<\beta<z'(T),$ then the next problem to solve is (\ref{numer}) with  $a = (y''(0)+z''(0))/2.$

{\bf Examples.}

While solving the initial value problems we apply an adaptive Runge-Kutta-Fehlberg method RK45 \cite{Mat}, in which a tolerance parameter $\epsilon$ controls local error of the method.  The values of $\epsilon$ range from $10^{-8}$ to $10^{-14}.$ For representing real values we use standard $16-17$-digits double data type.

{\bf Numerical results for the classical Blasius equation $p=1,$ $c=1/2,$ and $\beta=1.$}

Here $a_{min}=0.2694860459,$ $a_{max}=0.3420953216.$ For
$T=14$ we get the all three inequalities (\ref{2-der}), (\ref{1-der}) and (\ref{0-der}) for $\epsilon=1.0e-14.$
 Below $N$ stands for a number of steps in RK45 (average):
\vspace{.3cm}

\begin{tabular}{|l|l|l|l|l|l|}
  \hline
  $\epsilon$ & $N$ & $a$ & $|x''(0)-a|$ & $|1-x'(T)|$ & $x(T)$ \\
  \hline
  1.0e-08 & 112 & 0.332057330068201 & 6.15e -09 & 1.2e-08 & 12.279212180321 \\
  1.0e-09 & 197 & 0.332057335646357 & 5.69e -10 & 1.1e-09 & 12.279212327474 \\
  1.0e-10 & 338 & 0.332057336149237 & 6.60e -11 & 1.3e-10 & 12.279212340740 \\
  1.0e-11 & 611 & 0.332057336210248 & 4.95e -12 & 9.9e-12 & 12.279212342350 \\
  1.0e-12 & 971 & 0.332057336214903 & 2.94e -13 & 5.7e-13 & 12.279212342472 \\
  1.0e-13 & 1831 & 0.332057336215154 & 4.19e-14 & 6.5e-14 & 12.279212342479 \\
  1.0e-14 & 3346 & 0.332057336215186 & 1.06e-14 & 5.5e-16 & 12.279212342480 \\
  \hline
\end{tabular}

\vspace{.4cm}
The last value of $a$ differs in two last digits from the one cited in \cite{Bo1}: $$a=0.33205733621519630.$$
The last two columns of the table have been computed for $\epsilon=10^{-14}.$ As there is $x''(T) = 7.68e-13,$ for $t>T$ the straight line approximation of the solution is the most effective.

{\bf Numerical results for the equation with $p=7,$ $c=1/2$ and $\beta=1.$}

Here, $a_{min}=0.3733978388,$ $a_{max}=0.3805482427.$ As above for the tolerance $\epsilon=1.0e-14,$ the interval $[0,4]$ is sufficiently large and we get the following
results by RK45 method:
\vspace{.3cm}

\begin{tabular}{|l|l|l|l|l|l|}
  \hline
  $\epsilon$ & $N$ & $a$ & $|x''(0)-a|$ & $|1-x'(T)|$ & $x(T)$ \\
  \hline
  1.0e-08 & 189 & 0.379398164451122 & 2.47e-08 & 3.5e-08 & 2.673055448977 \\
  1.0e-09 & 316 & 0.379398187063634 & 2.04e-09 & 2.9e-09 & 2.673055570853 \\
  1.0e-10 & 549 & 0.379398189005442 & 1.03e-10 & 1.4e-10 & 2.673055581319 \\
  1.0e-11 & 961 & 0.379398189086642 & 2.20e-11 & 3.1e-11 & 2.673055581757 \\
  1.0e-12 & 1688 & 0.379398189106905 & 1.69e-12 & 2.3e-12 & 2.673055581866 \\
  1.0e-13 & 2827 & 0.379398189108438 & 1.62e-13 & 1.9e-13 & 2.673055581874 \\
  1.0e-14 & 4634 & 0.379398189108571 & 2.91e-14 & 0.0e-14 & 2.673055581875 \\
  \hline
\end{tabular}

\vspace{.4cm}

Remarks –- as above. Here $x''(T) =9.03e-18.$  The value of $a= 0.3793981891086$ -- here, all digits are true.

{\bf Numerical experiment for the equation with $p=0.1$ $c=1/2,$ $\beta=1;$} taking $\epsilon=10^{-14}$ we get $T=50.$ The proof of the existence and uniqueness result
for $p<1$ fails, since we cannot claim that the initial value problem (\ref{equ-ini}) has a unique solution and that it depends continuously on $a.$ Hence,
function $h$ can be multivalued. If one will prove the uniqueness, then, due to \cite{PSV} p. 172, $h$ will be continuous and all results of this paper will be true also for $p<1.$ The stability of numerical experiments cited below suggests it is the fact.

\vspace{.3cm}

\begin{tabular}{|l|l|l|l|l|l|}
  \hline
  $\epsilon$ & $N$ & $a$ & $|x''(0)-a|$ & $|1-x'(T)|$ & $x(T)$ \\
  \hline
  1.0e-08 & 142 & 0.443643205985427 & 2.16e-07 & 4.5e-07 & 48.05426086324 \\
  1.0e-09 & 256 & 0.443643403844908 & 1.78e-08 & 3.7e-08 & 48.05428058120 \\
  1.0e-10 & 466 & 0.443643420192529 & 1.49e-09 & 3.1e-09 & 48.05428221034 \\
  1.0e-11 & 839 & 0.443643421402885 & 2.80e-10 & 5.8e-10 & 48.05428233096 \\
  1.0e-12 & 1505 & 0.443643421660499 & 2.25e-11 & 4.7e-11 & 48.05428235664 \\
  1.0e-13 & 2669 & 0.443643421681506 & 1.49e-12 & 3.6e-12 & 48.05428235873 \\
  1.0e-14 & 4922 & 0.443643421683245 & 2.45e-13 & 2.2e-16 & 48.05428235890 \\
  \hline
\end{tabular}

\vspace{.4cm}

Remarks – as above. Here $x''(T) =1.02e-15$. The value of $a= 0.443643421683$ -- here, all digits are true.

\section{Conclusions}
The authors know that our computation of the value of the second derivative of the solution are not more exact than others.
However, the proposed method gives a possibility of controlling errors and it is very simple. We hope a similar approach
can be applied for more general equations as $x'''+f(x)\cdot g(x'')=0$ with qualitative assumptions on functions $f$ and $g.$

\noindent Address:
Institute of Mathematics\\
\hbox{\vspace{1cm}} Technical Univ. of Lodz\\
\hbox{\vspace{1cm}} Wolczanska 215, 90-924 Lodz,\\
\hbox{\vspace{1cm}} Poland,\\
e-mail: gandrzejczak99@wp.pl, magdan@p.lodz.pl, przeradz@p.lodz.pl.


\begin{thebibliography}{000}


\bibitem{BL} H. Blasius, {\em Grenzschichten in Fl\"ussigkeiten mit kleiner Reibung}, Math. Phys., 56 (1908), 1--37.

\bibitem{Bo1} J. P. Boyd, {\em The Blasius function in the complex plane}, Exp. Math., 8 (1999), 381--394.

\bibitem{Bo2} J. P. Boyd, {\em The Blasius function: computations before computers, the value of tricks, undergraduate projects and open research problems}, SIAM Review, 50 (2008), 791--804.

\bibitem{BH} B. Brighi, J.-D. Hoernel {\em On the concave and convex solutions of a mixed convection
boundary layer approximation in a porous medium}, Appl. Math. Lett. 19 (2006) 69–-74.

\bibitem{Mat} J. H. Mathews, {\em Numerical Methods for Computer Science, Engineering and Mathematics}, Prentice Hall Inc., New Jersey 1992.

\bibitem{PT} K. Parand, A. Taghavi, {\em Rational scaled generalized Laguerre function collocation method for solving the Blasius equation}, J.  Comp. and Appl. Math., 233 (2009), 980--989.

\bibitem{PSV} L. C. Piccinini, G. Stampacchia, G. Vidossich, {\em Ordinary Differential Equations in $\R ^n$}, Applied Mathematical Sciences 39, Springer-Verlag, New York-Berlin-Heidelberg 1984.

\bibitem{Ran} A. I. Ranasinghe, {\em Solution of Blasius equation by decomposition}, Appl. Math. Sci., 3 (2009), 605-611.

\end{thebibliography}
\end{document}